# NONPARAMETRIC ESTIMATION OVER SHRINKING NEIGHBORHOODS: SUPEREFFICIENCY AND ADAPTATION[1]

By T. Tony Cai and Mark G. Low

*University of Pennsylvania*

A theory of superefficiency and adaptation is developed under flexible performance measures which give a multiresolution view of risk and bridge the gap between pointwise and global estimation. This theory provides a useful benchmark for the evaluation of spatially adaptive estimators and shows that the possible degree of superefficiency for minimax rate optimal estimators critically depends on the size of the neighborhood over which the risk is measured.

Wavelet procedures are given which adapt rate optimally for given shrinking neighborhoods including the extreme cases of mean squared error at a point and mean integrated squared error over the whole interval. These adaptive procedures are based on a new wavelet block thresholding scheme which combines both the commonly used horizontal blocking of wavelet coefficients (at the same resolution level) and vertical blocking of coefficients (across different resolution levels).

**1. Introduction.** Squared error loss at each point and integrated squared error loss over an interval are two of the most common ways to evaluate the performance of nonparametric function estimators. Integrated squared error is used as a broad overall measure of loss whereas pointwise squared error loss gives a highly localized measure of accuracy. Minimax theory for both these cases can be found for example in Pinsker (1980), Ibragimov and Hasminski (1984), Donoho and Liu (1991) and Donoho and Johnstone (1998), and there are a large number of additional references in Efromovich (1999).

In nonparametric function estimation problems minimax risk provides a useful uniform benchmark for the comparison of estimators. Such uniform bounds do not, however, capture many aspects of these problems since in these infinite-dimensional settings asymptotically minimax estimators can

Received December 2002; revised March 2004.
[1]Supported in part by NSF Grant DMS-03-06576.

often be constructed which are also superefficient at every parameter point. In fact, much recent work on nonparametric function estimation can be viewed as attempts to construct superefficient estimators with desirable properties. This is clear in the literature on adaptive estimation where the connection between superefficiency and adaptation has been considered as, for example, in Beran (1999, 2000). In adaptive estimation the goal is to construct estimators which are simultaneously asymptotically (near) minimax over a collection of parameter spaces. Such estimators are optimal over this range of spaces.

This theory of adaptive estimation depends strongly on how risk is measured. When the performance is measured globally full adaptation can often be achieved. In particular, Efromovich and Pinsker (1984) constructed fully adaptive estimators over a range of Sobolev spaces. Recent results on rate adaptive estimators focus on more general Besov spaces. See, for example, Donoho and Johnstone (1995), Cai (1999) and Härdle, Kerkyacharian, Picard and Tsybakov (1998).

When the performance is measured at a point, it is often the case that full adaptation is not possible and superefficient estimators must have inflated risk at other parameter points. A penalty, usually a logarithmic factor, must be paid for not knowing the smoothness. Important work in this area began with Lepski (1990) where attention focused on a collection of Lipschitz classes. See also Brown and Low (1996), Efromovich and Low (1994), Lepski and Spokoiny (1997) and Cai (2003).

Since optimally adaptive estimators at each point typically pay a logarithmic penalty compared to the minimax risk, they are not necessarily optimally globally adaptive. This has led to the approach of a simultaneous pointwise and global analysis. The goal is then to construct estimators which, for a range of parameter spaces, are both minimax rate optimal for integrated squared error loss and pay only a logarithmic penalty for squared error loss at each point. See, for example, Cai (1999, 2002) and Efromovich (2002).

Pointwise mean squared error can be viewed as an extreme (although useful) way of measuring local performance of an estimator $\hat{f}_n$. The focus in the present paper is on a more flexible approach. Specifically we propose to evaluate the performance of an estimator $\hat{f}_n$ (near $x_0$) by using an average mean squared error over a neighborhood of $x_0$:

$$(1) \qquad R(\hat{f}_n, f; x_0, c_n) \equiv \frac{1}{2c_n} E_f \int_{x_0-c_n}^{x_0+c_n} (\hat{f}_n(x) - f(x))^2 \, dx.$$

The choice of $c_n$ allows for considerable flexibility when measuring local performance. For fixed $n$, by taking the limit as $c_n \to 0$ we can recover the usual case of squared error loss at $x_0$, and by taking $x_0 = \frac{1}{2}$ and $c_n = \frac{1}{2}$ we recover the usual global risk. By evaluating the performance for a whole range of $c_n$ it is possible to give a multiresolution view of the risk. We show that this more flexible approach to measuring local performance can be used to bridge the gap between the pointwise and global theories.

In this paper we consider estimation over shrinking neighborhoods based on observations from a Gaussian process

$$(2) \qquad Z_n^*(t) \equiv \int_0^t f(x)\,dx + \frac{1}{\sqrt{n}} B^*(t), \qquad 0 \le t \le 1,$$

where $B^*(t)$ is a standard Brownian motion and $f$ is an unknown function. This Gaussian process is a prototypical model for many nonparametric function estimation problems such as nonparametric regression and density estimation.

In Section 2 it is shown that the size of the neighborhood as governed by $c_n$ determines both the possible degree of superefficiency for minimax rate optimal estimators as well as the cost of adaptation. For "small" neighborhoods superefficient estimators cannot be minimax rate optimal and hence fully rate adaptive estimation is not possible. In fact the penalty for superefficiency determines the cost of adaptation. On the other hand, for "large" neighborhoods there exist minimax rate optimal estimators which are superefficient at every parameter point.

Adaptive estimation is considered in Sections 3 and 4. In Section 3 a procedure is constructed which optimally adapts to smoothness over given shrinking neighborhoods. This construction includes the extreme cases of mean squared error at a point and mean integrated squared error over the whole interval.

The adaptive procedure used in Section 3 is based on block thresholding of empirical wavelet coefficients, a technique which has been shown to be effective for adaptive estimation. See, for example, Hall, Kerkycharian and Picard (1998) and Cai (1999, 2002). Block thresholding in these papers is done by blocking of wavelet coefficients only at the same resolution level. The adaptive procedure proposed here is based on a new block thresholding scheme. It combines both the commonly used horizontal blocking of wavelet coefficients (at the same resolution level) and vertical blocking of coefficients (across different resolution levels). Furthermore, it appears that vertical blocking is essential for the resulting estimator to be optimally adaptive.

The theory of adaptive estimation over given shrinking neighborhoods developed in Sections 2 and 3 provides a useful benchmark for the evaluation of estimators designed to be spatially adaptive. Spatially adaptive procedures should however adapt not just to unknown smoothness but also to a whole range of shrinking neighborhoods over which the risk is measured. This more complete analysis incorporating a multi-resolution view of risk is given in Section 4. In that section it is shown that a block thresholding estimator introduced in Cai (1999) exhibits, from this point of view, good spatial adaptivity.

**2. Superefficiency and adaptation.** In nonparametric function estimation problems minimax risk depends strongly on the parameter space. Typically the parameter space is unknown and so attention is often focused on the construction of adaptive estimators which simultaneously attain near minimaxity over a collection of parameter spaces. The theory of adaptive estimators is closely connected to that of superefficient estimators which in turn depend on how the risk is measured.

In this paper we shall develop the shrinking neighborhood theory for Hölder classes

(3) $\quad F(\alpha, M) = \{f : |f^{(k)}(x) - f^{(k)}(y)| \leq M|x-y|^{\alpha-k}, 0 \leq x < y \leq 1\},$

where $k$ is the greatest integer strictly less than $\alpha$. Minimax theory in this setup is standard. In particular, under the risk measure (1) with observations from the Gaussian process (2) the minimax rate of convergence over $F(\alpha, M)$ is of order $n^{-2\alpha/(2\alpha+1)}$. The theory for superefficiency and adaptation is however quite interesting.

The focus in this section is on how the size of the shrinking neighborhood affects the penalty for superefficient estimators. The connection between superefficient estimation and adaptation is then made clear. Our interest in superefficiency is mainly for the insight it provides for the question of adaptation and we show how lower bounds derived for the penalty of superefficiency are directly applicable to the minimum cost of adaptation.

2.1. *Superefficiency.* For a parameter space $\mathcal{F}$ we call an estimator $\hat{f}_n$ superefficient at $f \in \mathcal{F}$ under a loss function $L_n$ if the risk at $f$ converges faster than the minimax risk, namely

(4) $$\frac{E_f L_n(\hat{f}_n, f)}{\inf_{\hat{f}_n} \sup_{f \in \mathcal{F}} E_f L_n(\hat{f}_n, f)} \to 0.$$

As mentioned in the Introduction, for estimation under mean integrated squared error (i.e., $x_0 = \frac{1}{2}$ and $c_n = \frac{1}{2}$) fully rate adaptive estimators exist and so there are superefficient estimators which are also minimax rate optimal. In particular, Brown, Low and Zhao (1997) give examples of estimating the whole function under integrated squared error loss where an estimator is superefficient at every parameter while also maintaining the minimax rate of convergence. On the other hand, for estimation under pointwise mean squared error ($c_n = 0$) Lepski (1990) and Brown and Low (1996) showed that any superefficient estimator cannot be minimax rate optimal over $F(\alpha, M)$ and hence in this case fully rate optimal adaptation is not possible. This case is similar to the superefficiency phenomenon arising in regular parametric problems. See, for example, Le Cam (1953) and Lehmann (1983).

As argued in the Introduction, integrated squared error and pointwise squared error are two extremes of a whole range of risk measures, each of which sheds light on the performance of a particular estimator $\hat{f}_n$. Shrinking neighborhoods give a more general way to evaluate the performance of an estimator. We begin by exploring the minimal cost of superefficiency for a specified shrinking neighborhood and find the critical size of neighborhood which will allow for the construction of superefficient estimators which are also minimax rate optimal.

For a given shrinking neighborhood of $x_0$ let $\Delta(f_0)$ be the collection of estimators $\hat{f}_n$ based on the Gaussian observations (2) that are superefficient at rate $B_n$ at the parameter point $f_0$. More specifically, let

$$(5) \qquad \Delta(f_0) = \left\{ \hat{f}_n : \limsup_{n \to \infty} n^{2\alpha/(1+2\alpha)} B_n R(\hat{f}_n, f_0; x_0, c_n) < \infty \right\}.$$

The following result then precisely quantifies the minimum penalty of such superefficient estimators.

THEOREM 1. *Fix $0 < x_0 < 1$, $0 < M' < M$, and set $c_n = d_n n^{-1/(1+2\alpha)}$. Let $B_n \to \infty$ and $\frac{n}{\log B_n} \to \infty$ and suppose that $f_0 \in F(\alpha, M')$.*

(i) *If $\limsup_{n \to \infty} d_n \cdot (\log B_n)^{-1/(1+2\alpha)} = 0$, then for any $\hat{f}_n \in \Delta(f_0)$,*

$$(6) \qquad \liminf_{n \to \infty} \left( \frac{n}{\log B_n} \right)^{2\alpha/(1+2\alpha)} \sup_{f \in F(\alpha, M)} R(\hat{f}_n, f; x_0, c_n) > 0$$

*and there exists some $\hat{f}_n \in \Delta(f_0)$ satisfying*

$$(7) \qquad \limsup_{n \to \infty} \left( \frac{n}{\log B_n} \right)^{2\alpha/(1+2\alpha)} \sup_{f \in F(\alpha, M)} R(\hat{f}_n, f; x_0, c_n) < \infty.$$

(ii) *If $\liminf_{n \to \infty} d_n \cdot (\log B_n)^{-1/(1+2\alpha)} > 0$ and $\limsup_{n \to \infty} d_n \times (\log B_n)^{-1} = 0$, then for any $\hat{f}_n \in \Delta(f_0)$,*

$$(8) \qquad \liminf_{n \to \infty} n^{2\alpha/(1+2\alpha)} \cdot \frac{d_n}{\log B_n} \sup_{f \in F(\alpha, M)} R(\hat{f}_n, f; x_0, c_n) > 0$$

*and there exists some $\hat{f}_n \in \Delta(f_0)$ satisfying*

$$(9) \qquad \limsup_{n \to \infty} n^{2\alpha/(1+2\alpha)} \frac{d_n}{\log B_n} \sup_{f \in F(\alpha, M)} R(\hat{f}_n, f; x_0, c_n) < \infty.$$

(iii) *If* $\liminf_{n\to\infty} \frac{d_n}{\log B_n} > 0$, *then there exists an estimator* $\hat{f}_n \in \Delta(f_0)$ *satisfying*

(10) $$\limsup_{n\to\infty} n^{2\alpha/(1+2\alpha)} \sup_{f\in F(\alpha,M)} R(\hat{f}_n, f; x_0, c_n) < \infty.$$

Note that the rate in the upper bound in case (iii) is sharp because it is also the minimax rate of convergence.

Theorem 1 gives bounds on the maximum risk after prespecifying the degree of superefficiency. For each of the three cases the proof of Theorem 1 constructs specific wavelet block thresholding procedures which attain the lower bounds. In other words, these wavelet procedures have minimal maximum risk given a particular level of superefficiency at a specified function.

Alternatively, it is also useful to classify the existence of minimax superefficient estimators in terms of a given neighborhood. The results can then be conveniently summarized as follows.

CASE 1 (Small neighborhoods). When the size of the neighborhood is smaller than $Dn^{-1/(2\alpha+1)}$ (i.e., $0 \le d_n \le D$) for some constant $D$, no minimax rate optimal estimator can be superefficient. In particular, when $d_n = 0$, which corresponds to the usual pointwise risk at $x_0$, superefficient estimators cannot be minimax rate optimal. In other words, minimax rate optimal estimators must have the same "flat" rate of convergence at every $f$ in the interior of $F(\alpha, M)$.

CASE 2 (Large neighborhoods). When the size of the neighborhood satisfies $\liminf d_n = \infty$ there are superefficient estimators attaining the minimax rate. The possible degree of superefficiency of a minimax rate optimal estimator however depends on the size of the neighborhood as described in the following three cases.

*Case* A. $\liminf_{n\to\infty} d_n = \infty$ and $\limsup_{n\to\infty} \frac{d_n}{\log n} = 0$. In this case a minimax rate optimal estimator can be superefficient at $f_0$, but the rate of convergence of its risk at $f_0$ cannot be algebraically faster than the minimax rate.

*Case* B. $0 < \liminf_{n\to\infty} \frac{d_n}{\log n} < \limsup_{n\to\infty} \frac{d_n}{\log n} < A < \infty$. In this case an estimator can have risk at $f_0$ converging at an algebraic rate faster than the minimax rate while maintaining the minimax convergence rate over $F(\alpha, M)$.

*Case* C.  $\liminf_{n\to\infty} \frac{d_n}{\log n} = \infty$. In this case a minimax rate optimal estimator can have its risk at $f_0$ converging at a rate which is faster than any algebraic rate. Hence an estimator can achieve a high degree of superefficiency at $f_0$ without paying a penalty in terms of its maximum risk over $F(\alpha, M)$.

An interesting consequence of these results is that for a prespecified shrinking neighborhood of size $n^{-\gamma}$ superefficient estimators which are also minimax rate optimal exist for $F(\alpha, M)$ if and only if $0 < \alpha < \frac{1-\gamma}{2\gamma}$. In particular, for $\gamma \geq 1$ there are no minimax superefficient estimators over any Hölder class $F(\alpha, M)$ and for $0 < \gamma < 1$ superefficient minimax rate optimal estimators exist only for the less smooth function spaces.

2.2. *Superefficiency in global estimation.* An interesting special case of the results considered in the previous section is that of estimation under mean integrated squared error which corresponds to the choice of $x_0 = \frac{1}{2}$ and $c_n = \frac{1}{2}$. In this case the results of Theorem 1 show that an estimator can simultaneously attain the minimax rate over $F(\alpha, M)$ and a high degree of superefficiency at any specific $f_0$ in the interior of $F(\alpha, M)$. The following corollary of Theorem 1 precisely quantifies how superefficient the estimator can be while maintaining the minimax rate of convergence over $F(\alpha, M)$.

COROLLARY 1.  *Let $0 < M' < M$ and $f_0 \in F(\alpha, M')$. Suppose*

$$\limsup_{n\to\infty} n^{1/(1+2\alpha)} \cdot (\log B_n)^{-1} = 0. \tag{11}$$

*If $\hat{f}_n$ is an estimator based on (2) satisfying*

$$\limsup_{n\to\infty} B_n E_{f_0} \|\hat{f}_n - f_0\|_2^2 < \infty, \tag{12}$$

*then*

$$\limsup_{n\to\infty} n^{2\alpha/(1+2\alpha)} \sup_{f \in F(\alpha, M)} E_f \|\hat{f}_n - f\|_2^2 = \infty. \tag{13}$$

Thus, a minimax rate optimal estimator cannot have risk at $f_0$ converging faster than $e^{-Dn^{1/(1+2\alpha)}}$ for all $D > 0$.

Condition (11) is sharp. That is, there exist estimators which converge super-fast at any fixed $f_0 \in F(\alpha, M)$ with the rate of $e^{-Dn^{1/(1+2\alpha)}}$ and yet still attain the minimax rate uniformly over the class $F(\alpha, M)$.

THEOREM 2.  *Let $f_0 \in F(\alpha, M)$ be fixed. For any constant $D > 0$ there exists an estimator which satisfies*

$$\limsup_{n\to\infty} e^{Dn^{1/(1+2\alpha)}} E_{f_0} \|\hat{f}_n - f_0\|_2^2 < \infty \tag{14}$$

*and*

(15) $$\limsup_{n\to\infty} n^{2\alpha/(1+2\alpha)} \sup_{f\in F(\alpha,M)} E_f \|\hat{f}_n - f\|_2^2 < \infty.$$

The theorem guarantees the existence of such superefficient estimators. One such estimator based on block thresholding of empirical wavelet coefficients is given by (63) and (64) in Section 6.

2.3. *Connection to adaptation.* The results on superefficiency given in Section 2.1 have direct implications for adaptation. Consider two function classes $F(\alpha_1, M)$ and $F(\alpha_2, M)$ with $0 < \alpha_1 < \alpha_2 \leq 1$. Then $F(\alpha_2, M) \subset F(\alpha_1, M)$ and a fully rate adaptive estimator $\hat{f}_n$ over these classes would need to satisfy

(16) $$\sup_{f\in F(\alpha_i, M)} R(\hat{f}_n, f; x_0, c_n) \asymp n^{-2\alpha_i/(2\alpha_i+1)}$$

for both $i=1$ and $i=2$. The risk of $\hat{f}_n$ for each $f \in F(\alpha_2, M)$ must then converge faster than the minimax risk over the larger parameter space $F(\alpha_1, M)$. Hence such estimators must be superefficient at each $f \in F(\alpha_2, M)$ with respect to $F(\alpha_1, M)$. The results in Theorem 1 can then be applied to yield corresponding lower bounds for adaptation over shrinking neighborhoods. These results are summarized in the following corollary.

COROLLARY 2. *Consider two function classes $F(\alpha_1, M_1)$ and $F(\alpha_2, M_2)$ with $\alpha_1 < \alpha_2$. Let $0 < x_0 < 1$ and $c_n = d_n n^{-1/(1+2\alpha_1)}$. If $\limsup_{n\to\infty} d_n \times (\log n)^{-1} = 0$, then*

(17) $$\max_{i=1,2} \limsup_{n\to\infty} n^{2\alpha_i/(1+2\alpha_i)} \inf_{\hat{f}_n} \sup_{f\in F(\alpha_i, M_i)} R(\hat{f}_n, f; x_0, c_n) = \infty.$$

*More specifically, suppose $\hat{f}_n$ is any estimator satisfying*

(18) $$\limsup_{n\to\infty} n^r \sup_{f\in F(\alpha_2, M_2)} R(\hat{f}_n, f; x_0, c_n) < \infty$$

*for some $r > \frac{2\alpha_1}{1+2\alpha_1}$.*

(i) *If $\limsup_{n\to\infty} d_n \cdot (\log n)^{-1/(1+2\alpha_1)} = 0$, then*

(19) $$\liminf_{n\to\infty} \left(\frac{n}{\log n}\right)^{2\alpha_1/(1+2\alpha_1)} \sup_{f\in F(\alpha_1, M_1)} R(\hat{f}_n, f; x_0, c_n) > 0.$$

(ii) *If $\liminf_{n\to\infty} d_n \cdot (\log n)^{-1/(1+2\alpha_1)} > 0$ and $\limsup_{n\to\infty} d_n \cdot (\log n)^{-1} = 0$, then*

(20) $$\liminf_{n\to\infty} n^{2\alpha_1/(1+2\alpha_1)} \frac{d_n}{\log n} \sup_{f\in F(\alpha_1, M_1)} R(\hat{f}_n, f; x_0, c_n) > 0.$$

The results in this corollary state that it is impossible to adaptively attain the minimax rates over the two function classes with different convergence rates whenever the size of the neighborhood is "too small." In Section 3 it is shown that the lower bounds on the cost of adaptation given by (19) and (20) are in fact sharp.

**3. Adaptive estimation.** We now turn our attention to adaptive estimation and the construction of adaptive estimators. In this section the focus is on adaptation over smoothness classes for a given shrinking neighborhood. Wavelet thresholding estimators are constructed which attain the bounds given by (19) and (20). In Section 4 we shall consider adaptation to both smoothness and to the size of the neighborhood.

3.1. *Wavelet thresholding.* Let $\phi$ and $\psi$ be a pair of compactly supported father and mother wavelets which generate an orthonormal basis of $L^2[0,1]$ through dilation and translation and where as is typical $\phi$ is chosen to satisfy $\int \phi = 1$. The support lengths of $\phi$ and and $\psi$ are written as $N_\phi$ and $N_\psi$, respectively.

Throughout this paper it is also assumed that $\psi$ is $r$-regular, meaning it has $r$ vanishing moments and $r$ continuous derivatives. Under these assumptions let

$$\phi_{j,k}(t) = 2^{j/2}\phi(2^j t - k), \qquad \psi_{j,k}(t) = 2^{j/2}\psi(2^j t - k).$$

Then for some choice of $j_0$ the collection $\{\phi_{j_0,k}, k = 1, \ldots, 2^{j_0}; \psi_{j,k}, j \geq j_0, k = 1, \ldots, 2^j\}$ with appropriate boundary corrections is an orthonormal basis of $L^2[0,1]$. See Cohen, Daubechies, Jawerth and Vial (1993), Daubechies (1994) and Meyer (1991) for further details on wavelet bases on the unit interval $[0,1]$. For wavelets on the line, see Daubechies (1992) and Meyer (1992).

A function $f:[0,1] \to \mathbb{R}$ can then be expanded in this orthonormal series. Set

$$\xi_{j_0,k} = \int_0^1 f(t)\phi_{j_0,k}(t)\,dt \quad \text{and} \quad \theta_{j,k} = \int_0^1 f(t)\psi_{j,k}(t)\,dt,$$

where $\xi_{j_0,k}$ are the wavelet coefficients at the coarse level and $\theta_{j,k}$ are coefficients at the detail levels.

Under the orthonormal wavelet basis, the Gaussian model (2) is equivalent to the sequence model

(21) $\qquad \tilde{y}_{j_0,k} = \xi_{j_0,k} + n^{-1/2}\tilde{z}_{j_0,k}, \qquad 1 \leq k \leq 2^{j_0},$

(22) $\qquad y_{j,k} = \theta_{j,k} + n^{-1/2}z_{j,k}, \qquad 1 \leq k \leq 2^j, j \geq j_0,$

where $\tilde{z}_{j_0,k}$ and $z_{j,k}$ are i.i.d. $N(0,1)$ random variables. Reconstructions of $f$ can then be based on estimates of the wavelet coefficients.

One particularly effective technique for estimating the wavelet coefficients is that based on block thresholding. Block thresholding estimates the wavelet coefficients in groups rather than individually, making simultaneous decisions to retain or to discard all the empirical coefficients within a block. It increases estimation accuracy by using information about neighboring wavelet coefficients balancing variance and bias along the curve. More details of such adaptive smoothing can be found in Hall, Kerkyacharian and Picard (1998) and Cai (1999, 2002). More standard term-by-term thresholding rules can be thought of as a special case of block thresholding with block size one.

The block thresholding rules used in the above-mentioned papers are constructed by grouping wavelet coefficients only at the same resolution level. In our context it is necessary to use block thresholding rules which employ vertical blocking of coefficients across different resolution levels as well as the commonly used horizontal blocking of wavelet coefficients at the same resolution level. We thus give below a generic description of a general block thresholding estimator which possibly uses both horizontal and vertical blocking.

Let $J > j_0$ be some dividing resolution level. Group the wavelet coefficients from level $j_0$ to level $J$ into nonoverlapping blocks of length $L$. Let $B_i$ be the set of indices for coefficients in the $i$th block and let $S_i^2 = \sum_{(j,k) \in B_i} y_{j,k}^2$ be the sum of squares for this block. The block thresholding estimator of the wavelet coefficients has the form

$$\hat{\theta}_{j,k} = \begin{cases} \eta(S_i^2) \cdot y_{j,k}, & \text{for } (j,k) \in B_i, j \leq J, \\ 0, & \text{for } j > J, \end{cases} \tag{23}$$

where $\eta(S_i^2)$ is some thresholding function. For example, one can take $\eta(S_i^2) = I(S_i^2 > \lambda)$ or

$$\eta(S_i^2) = \left(1 - \frac{\lambda L n^{-1}}{S_i^2}\right)_+ \tag{24}$$

where $\lambda$ is some thresholding constant. The shrinkage rule (24) is used throughout this paper with a variety of values of $\lambda$ and $L$.

3.2. *Adaptive estimation on given neighborhoods.* The lower bound on the performance of an adaptive estimator over a collection of Hölder classes $\mathcal{F}(\alpha, M)$ has been given in Corollary 2 of Section 2.3. For neighborhoods with $c_n \leq \frac{\log n}{n}$ an estimator is given in Section 4 which adapts both to smoothness and to the size of neighborhood while attaining the bounds of Corollary 2. In this section, an estimator designed for neighborhoods with $c_n > \frac{\log n}{n}$ is given. It is a wavelet estimator based on a block thresholding scheme. Using the same notation as in Section 3.1, let $J$, $J_*$ and $J^*$ be the smallest integers satisfying

$$2^J \geq n, \qquad 2^{J_*} \geq c_n^{-1} \quad \text{and} \quad 2^{J^*} \geq c_n^{-1} \log n,$$

respectively. Then in the case $c_n > \frac{\log n}{n}$ considered here, it follows that $J^* < J$. We set $J_* = j_0$ when $J_* < j_0$ and let

$$H_j = \{(j,k) : \operatorname{supp}(\psi_{j,k}) \cap [x_0 - c_n, x_0 + c_n] \neq \varnothing\} \quad \text{and} \quad H^* = \bigcup_{J_* \leq j < J^*} H_j.$$

Then

$$\operatorname{Card}(H_j) \leq \begin{cases} N_\psi, & \text{if } j < J_*, \\ N_\psi 2^j c_n, & \text{if } j \geq J_*, \end{cases}$$

and $\operatorname{Card}(H^*) \asymp \log n$ where $N_\psi$ is the length of the support of $\psi$.

The estimator we propose, a hybrid estimator of soft thresholding, vertical block thresholding and horizontal block thresholding, can be described in four steps as follows.

1. For empirical coefficients $y_{j,k}$ between levels $j_0$ and $J_*$ apply term-by-term soft thresholding rule. The soft thresholding rule is also applied to coefficients at levels between $J_*$ and $J^*$ where $(j,k) \notin H^*$, in which case the support of the corresponding wavelet basis function $\psi_{j,k}$ has empty intersection with the interval $[x_0 - c_n, x_0 + c_n]$.
2. Group all the empirical coefficients $y_{j,k}$ with $(j,k) \in H^*$ into a single vertical block and denote by $S_v^2 = \sum_{(j,k) \in H^*} y_{j,k}^2$ the sum of squared coefficients in the vertical block. Apply a single James–Stein shrinkage rule of the form (24) to the coefficients in this block.
3. At each resolution level $J^* \leq j < J$, divide the empirical wavelet coefficients $y_{j,k}$ into nonoverlapping blocks of length $L = \log n$. Denote by $(jb)$ the set of indices of the coefficients in the $b$th block at level $j$, that is, $(jb) = \{(k : (b-1)L + 1 \leq k \leq bL\}$, and let $S_{(jb)}^2 = \sum_{k \in (jb)} y_{j,k}^2$ denote the sum of squares for the block $(jb)$. Then apply the James–Stein shrinkage rule to each block $(jb)$ for $J^* \leq j < J$.
4. For $j \geq J$, estimate all $\theta_{j,k}$ by 0.

More precisely, each coefficient $\theta_{j,k}$ is estimated by

$$(25) \quad \hat{\theta}_{j,k} = \begin{cases} \operatorname{sgn}(y_{j,k})(|y_{j,k}| - \sqrt{2n^{-1}\log n})_+, \\ \qquad\qquad \text{if } j_0 \leq j < J^* \text{ and } (j,k) \notin H^*, \\ \left(1 - \frac{\lambda_* L n^{-1}}{S_v^2}\right)_+ y_{j,k}, \quad \text{if } (j,k) \in H^*, \\ \left(1 - \frac{\lambda_* L n^{-1}}{S_{(jb)}^2}\right)_+ y_{j,k}, \quad \text{if } J^* \leq j < J \text{ and } k \in (jb), \\ 0, \qquad\qquad \text{if } j \geq J, \end{cases}$$

where $\lambda_* = 4.50524$ is a constant satisfying $\lambda_* - \log \lambda_* - 1 = 2$. Define the wavelet estimator of $f$ by

$$(26) \quad \hat{f}_n(x) = \sum_{k=1}^{2^{j_0}} \tilde{y}_{j_0,k} \phi_{j_0,k}(x) + \sum_{j=j_0}^{\infty} \sum_{k=1}^{2^j} \hat{\theta}_{j,k} \psi_{j,k}(x)$$

with $\hat{\theta}_{j,k}$ given in (25). This estimator attains the lower bounds in Corollary 2 at least when $c_n > \frac{\log n}{n}$. For smaller neighborhoods, the estimator given by (32) and (33) in Section 4 also attains the lower bounds. These results are summarized in the following theorem.

THEOREM 3. *When $c_n > \frac{\log n}{n}$, let $\hat{f}_n$ be the estimator given by (25) and (26) where the wavelet $\psi$ is $r$-regular with $r > \alpha$, whereas if $c_n \leq \frac{\log n}{n}$ let $\hat{f}_n$ be the BlockJS estimator given in (32) and (33). Let $0 < x_0 < 1$ and $c_n = d_n n^{-1/(1+2\alpha)}$.*

(i) *If $\limsup_{n\to\infty} d_n \cdot (\log n)^{-1/(1+2\alpha)} < \infty$, then*

$$\text{(27)} \qquad \limsup_{n\to\infty} \left(\frac{n}{\log n}\right)^{2\alpha/(1+2\alpha)} \sup_{f \in F(\alpha,M)} R(\hat{f}_n, f; x_0, c_n) < \infty.$$

(ii) *If $\liminf_{n\to\infty} d_n \cdot (\log n)^{-1/(1+2\alpha)} = \infty$ and $\limsup_{n\to\infty} d_n \times (\log n)^{-1} = 0$, then*

$$\text{(28)} \qquad \limsup_{n\to\infty} n^{2\alpha/(1+2\alpha)} \frac{d_n}{\log n} \sup_{f \in F(\alpha,M)} R(\hat{f}_n, f; x_0, c_n) < \infty.$$

(iii) *If $\liminf_{n\to\infty} d_n \cdot (\log n)^{-1} > 0$, then*

$$\text{(29)} \qquad \limsup_{n\to\infty} n^{2\alpha/(1+2\alpha)} \sup_{f \in F(\alpha,M)} R(\hat{f}_n, f; x_0, c_n) < \infty.$$

In view of Theorem 1, the estimator given in (25) and (26) attains the adaptive minimax rate for estimating $f$ over the neighborhood $[x_0 - c_n, x_0 + c_n]$.

A particularly interesting choice of $c_n$, $c_n = n^{-\gamma}$, is summarized in the following corollary which shows that fully rate optimal adaptation can be achieved over $F(\alpha, M)$ if and only if $0 < \alpha < \frac{1-\gamma}{2\gamma}$.

COROLLARY 3. *Let $\hat{f}_n$ be the estimator given in (25) and (26) and let $c_n = n^{-\gamma}$ for some $0 < \gamma < 1$. Suppose the wavelet $\psi$ is chosen to be $r$-regular with $r > \frac{1-\gamma}{2\gamma}$. Then for $0 < \alpha < \frac{1-\gamma}{2\gamma}$,*

$$\text{(30)} \qquad \limsup_{n\to\infty} n^{2\alpha/(1+2\alpha)} \sup_{f \in F(\alpha,M)} R(\hat{f}_n, f; x_0, c_n) < \infty,$$

*and for $\frac{1-\gamma}{2\gamma} \leq \alpha < r$,*

$$\text{(31)} \qquad \limsup_{n\to\infty} \left(\frac{n}{\log n}\right)^{2\alpha/(1+2\alpha)} \sup_{f \in F(\alpha,M)} R(\hat{f}_n, f; x_0, c_n) < \infty.$$

**4. Adaptation over smoothness and neighborhoods.** In nonparametric function estimation, it is common to fix a risk measure such as integrated squared error or squared error at a given point and to construct estimators which adapt across a range of smoothness classes. In our setting of shrinking neighborhoods, it is natural to consider two different types of adaptation. One is to adapt to the unknown smoothness of the underlying functions while the risk is measured over a given sequence of shrinking neighborhoods as in Section 3. A more ambitious and general adaptation goal is to adapt both to the unknown smoothness and the shrinking neighborhood over which the risk is measured.

This latter approach is most appropriate when the goal is to construct spatially adaptive estimators. It gives a more complete analysis with a multiresolution view of risk which spans a whole range of local and global measures of risk. Ideally we would like to construct an estimator which is "fully" adaptive—attaining the best adaptive rates for all choices of neighborhood sizes. The benchmark for such estimators is provided in Theorem 1. We shall show below that the BlockJS estimator [Cai (1999)] is nearly fully adaptive. This BlockJS procedure can be described as follows.

Expand the Gaussian process (2) in an orthonormal wavelet basis as in Section 3.1. At each resolution level $j < J = [\log_2 n]$ divide the empirical wavelet coefficients $y_{j,k}$ into nonoverlapping blocks of length $L = \log n$. Denote by $(jb)$ the set of indices of the coefficients in the $b$th block at level $j$, that is,

$$(jb) = \{k : (b-1)L + 1 \leq k \leq bL\}.$$

Let $S^2_{(jb)} = \sum_{k \in (jb)} y^2_{j,k}$ denote the sum of squares for the block $(jb)$ and let $\lambda_* = 4.50524$ be given as in Section 3, the root of the equation $\lambda - \log \lambda - 1 = 2$. We then apply the James–Stein shrinkage rule to each block $(jb)$ for $j_0 \leq j < J$,

$$(32) \qquad \hat{\theta}_{j,k} = \begin{cases} \left(1 - \dfrac{\lambda_* L n^{-1}}{S^2_{(jb)}}\right)_+ y_{j,k}, & \text{for } k \in (jb), j < J, \\ 0, & \text{for } j \geq J. \end{cases}$$

The BlockJS estimator $\hat{f}_n$ of the whole function $f$ is then given by

$$(33) \qquad \hat{f}_n(x) = \sum_{k=1}^{2^{j_0}} \tilde{y}_{j_0,k} \phi_{j_0,k}(x) + \sum_{j=j_0}^{\infty} \sum_{k=1}^{2^j} \hat{\theta}_{j,k} \psi_{j,k}(x).$$

THEOREM 4. *Let $\hat{f}_n$ be the BlockJS estimator given in* (32) *and* (33) *and let $0 < x_0 < 1$ and $c_n = d_n n^{-1/(1+2\alpha)}$. Suppose the wavelet $\psi$ is $r$-regular with $r > \alpha$.*

(i) If $\limsup_{n\to\infty} d_n \cdot (\log n)^{-1/(1+2\alpha)} < \infty$, then

$$\limsup_{n\to\infty} \left(\frac{n}{\log n}\right)^{2\alpha/(1+2\alpha)} \sup_{f\in F(\alpha,M)} R(\hat{f}_n, f; x_0, c_n) < \infty. \tag{34}$$

(ii) If $\liminf_{n\to\infty} d_n \cdot (\log n)^{-1/(1+2\alpha)} = \infty$ and $\limsup_{n\to\infty} d_n \cdot [(\log n) \times (\log\log n)]^{-1} = 0$, then

$$\limsup_{n\to\infty} n^{2\alpha/(1+2\alpha)} \frac{d_n}{(\log n)(\log\log n)} \sup_{f\in F(\alpha,M)} R(\hat{f}_n, f; x_0, c_n) < \infty. \tag{35}$$

(iii) If $\liminf_{n\to\infty} d_n \cdot [(\log n)(\log\log n)]^{-1} > 0$, then

$$\limsup_{n\to\infty} n^{2\alpha/(1+2\alpha)} \sup_{f\in F(\alpha,M)} R(\hat{f}_n, f; x_0, c_n) < \infty. \tag{36}$$

This theorem shows that the BlockJS estimator adapts well to the unknown smoothness across a wide range of shrinking neighborhoods. Just as in Section 3 the special choice $c_n = n^{-\gamma}$ is particularly interesting. Although the results of the following corollary are similar to those given in Corollary 3 it should be noted that the BlockJS estimator does not depend on the size or location of the neighborhood. Hence the BlockJS estimator exhibits very strong spatial and parameter space adaptivity.

COROLLARY 4. *Let $\hat{f}_n$ be the BlockJS estimator and let $0 < x_0 < 1$ and $c_n = n^{-\gamma}$ for some $\gamma > 0$. Suppose the wavelet $\psi$ is $r$-regular with $r \geq \frac{1-\gamma}{2\gamma}$. Then for $0 < \alpha < \frac{1-\gamma}{2\gamma}$,*

$$\limsup_{n\to\infty} n^{2\alpha/(1+2\alpha)} \sup_{f\in F(\alpha,M)} R(\hat{f}_n, f; x_0, c_n) < \infty, \tag{37}$$

*and for $\frac{1-\gamma}{2\gamma} \leq \alpha < r$,*

$$\limsup_{n\to\infty} \left(\frac{n}{\log n}\right)^{2\alpha/(1+2\alpha)} \sup_{f\in F(\alpha,M)} R(\hat{f}_n, f; x_0, c_n) < \infty. \tag{38}$$

**5. Discussion.** The theory of shrinking neighborhoods gives a multiresolution view of the performance of function estimators. It also provides a useful benchmark for the evaluation of spatially adaptive procedures. This theory can be easily extended to more general settings. One possible extension is to consider general weight functions. Let $w(x) \geq 0$ be a compactly supported continuous function satisfying $w(0) > 0$ and $\int w(x)\,dx = 1$. For a decreasing sequence $c_n \to 0$ and a fixed $x_0 \in (0,1)$ let

$$W_n(x) = \frac{1}{c_n} w\left(\frac{x - x_0}{c_n}\right). \tag{39}$$

The performance of an estimator $\hat{f}_n$ can then be evaluated with respect to the weight $W_n$:

$$(40) \qquad R(\hat{f}_n, f; W_n) \equiv E_f \int W_n(x)(\hat{f}_n(x) - f(x))^2 \, dx.$$

This risk can be viewed as a weighted risk concentrated around the point $x_0$, and the shrinking neighborhoods considered earlier in this paper correspond to the choice of uniform weight $w(x) = \frac{1}{2}I(-1 \leq x \leq 1)$.

Under the conditions given above, $w(x) \leq C_1$ for all $x$, $w(x) \geq C_2 > 0$ for $|x| \leq a$ and $w(x) = 0$ for $|x| \geq b$ for some constants $C_1$, $C_2$, $a$ and $b$. It is then easy to see that all the results given in the previous sections carry over to the risk given in (40).

It is also possible to extend the theory in this paper to a Gaussian process observed on the whole line. In this setting it is natural to consider a general weight function $W_n$ where $c_n \to 0$ or $c_n \to \infty$. The latter choice corresponds to expanding neighborhoods. When $c_n \to 0$ it is easy to see that all the theory given in the previous sections carries over to this setting. On the other hand, when $c_n \to \infty$ fully adaptive estimation is always possible and the block thresholding wavelet estimator given in Section 4 can easily be extended to a wavelet expansion on the real line.

**6. Proofs.** Throughout this section, $C$ denotes a generic positive constant which may vary from place to place. The wavelet notation follows that given in Section 3.1. The father wavelet $\phi$ and mother wavelet $\psi$ are always assumed to have compact support with the length of the support denoted by $N_\phi$ and $N_\psi$, respectively.

6.1. *Preparatory results.* The following elementary inequalities are useful for the evaluation of the risk of wavelet estimators over shrinking neighborhoods in terms of the wavelet coefficients.

LEMMA 1. *For any $0 \leq a < b \leq 1$, set*

$$S_1(a,b) \equiv \{(j,k) : \operatorname{supp}(\psi_{j,k}) \subset [a,b]\} \quad and$$

$$S_2(a,b) \equiv \{(j,k) : \operatorname{supp}(\psi_{j,k}) \cap [a,b] \neq \varnothing\}.$$

*Then*

$$(41) \qquad \sum_{(j,k) \in S_1(a,b)} \theta_{j,k}^2 \leq \int_a^b \left( \sum_{j,k} \theta_{j,k} \psi_{j,k}(x) \right)^2 dx \leq \sum_{(j,k) \in S_2(a,b)} \theta_{j,k}^2.$$

PROOF. For any $f(x) = \sum_{j,k} \theta_{j,k} \psi_{j,k}(x)$ let $h(x) = f(x) I_{[a,b]}(x)$ and note that $\int_a^b f^2(x) \, dx = \int_0^1 h^2(x) \, dx$. Let

$$g_1(x) = \sum_{(j,k) \in S_1(a,b)} \theta_{j,k} \psi_{j,k}(x) \quad \text{and} \quad g_2(x) = \sum_{(j,k) \in S_2(a,b)} \theta_{j,k} \psi_{j,k}(x).$$

Then $\|g_i\|_2^2 = \sum_{(j,k) \in S_i(a,b)} \theta_{j,k}^2$ for $i = 1, 2$. It is also easy to see that $g_2(x) = h(x)$ for $x \in [a,b]$ and so $\|g_2\|_2^2 \geq \|h\|_2^2$ and the second inequality in (41) immediately follows.

We can also write
$$h(x) = \sum_{j,k} \theta_{j,k} \psi_{j,k}(x) I_{[a,b]}(x) = g_1(x) + \sum_{(j,k) \notin S_1(a,b)} \theta_{j,k} \psi_{j,k}(x) I_{[a,b]}(x).$$

Noting that $\mathrm{supp}(g_1) \subset [a,b]$, it follows that
$$\int_0^1 g_1(x) \sum_{(j,k) \notin S_1(a,b)} \theta_{j,k} \psi_{j,k}(x) I_{[a,b]}(x) \, dx$$
$$= \int_0^1 g_1(x) \sum_{(j,k) \notin S_1(a,b)} \theta_{j,k} \psi_{j,k}(x) \, dx = 0,$$

and consequently
$$\|h\|_2^2 = \|g_1\|_2^2 + \left\| \sum_{(j,k) \notin S_1(a,b)} \theta_{j,k} \psi_{j,k}(x) I_{[a,b]}(x) \right\|_2^2 \geq \|g_1\|_2^2 = \sum_{(j,k) \in S_1(a,b)} \theta_{j,k}^2$$

and the first inequality in (41) also holds. □

The proofs of the main theorems also rely on bounds on the risk of wavelet block thresholding estimators. Lemma 2 summarizes several useful risk upper bounds for such estimators.

LEMMA 2. *Let $y_i = \theta_i + \sigma z_i$ where $z_i \stackrel{i.i.d.}{\sim} N(0,1), i = 1, \ldots, L$, and let $\hat{\theta}_i = (1 - \frac{\lambda L \sigma^2}{S^2})_+ y_i$ where $S^2 = \sum y_i^2$ and $\lambda \geq 1$. Then*

$$(42) \quad \sum_{i=1}^L E(\hat{\theta}_i - \theta_i)^2 \leq \min\left\{ \sum_{i=1}^L \theta_i^2, \lambda L \sigma^2 \right\} + 2\lambda e^{(-1/2)(\lambda - \log \lambda - 1)L} \sigma^2.$$

*In the special case of $\lambda = 4.50524$ (the root of the equation $\lambda - \log \lambda - 3 = 0$),*

$$(43) \quad \sum_{i=1}^L E(\hat{\theta}_i - \theta_i)^2 \leq \min\left\{ \sum_{i=1}^L \theta_i^2, \lambda L \sigma^2 \right\} + 2\lambda e^{-L} \sigma^2.$$

*In addition, suppose $\lambda = 4.50524$ and $|\theta_i| \leq c$ for all $i$. Then*

$$(44) \quad E(\hat{\theta}_i - \theta_i)^2 \leq 8c^2 + 2\lambda e^{-L} \sigma^2.$$

PROOF. Inequality (42) is a direct consequence of the oracle inequality given in Theorem 1 of Cai (1999) and the bound on the tail probability of the chi-squared distribution given in Lemma 2 of Cai (1999). Inequality (43) then follows directly on evaluation of (42). For the proof of (44), it suffices to consider the case of $\sigma = 1$. In that case note that

$$E(\hat{\theta}_i - \theta_i)^2 = E\left\{\left(1 - \frac{\lambda L}{S^2}\right)y_i I(S^2 > \lambda L) - \theta_i\right\}^2$$

$$\leq 2\theta_i^2 + 2Ey_i^2\left(1 - \frac{\lambda L}{S^2}\right)^2 I(S^2 > \lambda L)$$

$$\leq 2\theta_i^2 + 2Ey_i^2 I(S^2 > \lambda L).$$

Now note that for fixed $\theta_i$ it is easy to check that $Ey_i^2 I(S^2 > \lambda L)$ is increasing in each $|\theta_j|$ for $j \neq i$. Note also that if all $\theta_k$ other than $\theta_i$ are fixed, then by Lemma 3

(45) $$Ey_i^2 I(S^2 > \lambda L) = E(y_i^2 I(y_i^2 > \lambda L - S_{-i}^2)|S_{-i}^2)$$

is increasing in $|\theta_i|$. Hence $Ey_i^2 I(S^2 > \lambda L)$ is maximized when all $\theta_j = c$. When all $\theta_j = c$, $Ey_i^2 I(S^2 > \lambda L)$ is the same for all $i$ and hence

$$y_i^2 I(S^2 > \lambda L) = L^{-1} ES^2 I(S^2 > \lambda L).$$

In the proof of Proposition 1 in Cai (1999) it is shown that

$$ES^2 I(S^2 > \lambda L) \leq 3\|\theta\|_2^2 + \lambda L(\lambda^{-1}e^{\lambda-1})^{-L/2}.$$

Therefore,

$$E_\theta(\hat{\theta}_i - \theta_i)^2 \leq 8c^2 + 2\lambda e^{-(1/2)(\lambda - \log \lambda - 1)L}. \qquad \square$$

LEMMA 3. *Let $z \sim N(0,1)$ and $y = \theta + z$. Then for any $c \geq 0$, $E_\theta y^2 I(|y| > c)$ is an increasing function of $|\theta|$.*

PROOF. It suffices to consider $\theta \geq 0$. Let

$$f(\theta) = \sqrt{2\pi} E_\theta y^2 I(|y| > c) = \left(\int_c^\infty + \int_{-\infty}^{-c}\right) y^2 e^{-(1/2)(y-\theta)^2} \, dy$$

and

$$g(\theta) = \sqrt{2\pi} E_\theta z^2 I(|y| > c) = \left(\int_{c-\theta}^\infty + \int_{-\infty}^{-c-\theta}\right) z^2 e^{-(1/2)z^2} \, dz.$$

Then

$$f'(\theta) = \left(\int_c^\infty + \int_{-\infty}^{-c}\right) y^2 (y - \theta) e^{-(1/2)(y-\theta)^2} \, dy$$

$$= \left(\int_{c-\theta}^\infty + \int_{-\infty}^{-c-\theta}\right)(x^3 + 2\theta x^2 + \theta^2 x) e^{-(1/2)x^2} \, dx$$

and
$$g'(\theta) = (c-\theta)^2 e^{-(1/2)(c-\theta)^2} - (c+\theta)^2 e^{-(1/2)(c+\theta)^2}.$$

Note that
$$\int_{-\infty}^{-c-\theta} (x^3 + \theta^2 x) e^{-(1/2)x^2}\, dx = -\int_{c+\theta}^{\infty} (x^3 + \theta^2 x) e^{-(1/2)x^2}\, dx,$$

so for $\theta \geq 0$
$$f'(\theta) = \left(\int_{c-\theta}^{\infty} + \int_{-\infty}^{-c-\theta}\right) 2\theta x^2 e^{-(1/2)x^2}\, dx + \int_{c-\theta}^{c+\theta} (x^3 + \theta^2 x) e^{-(1/2)x^2}\, dx \geq 0,$$

and so the lemma follows. $\square$

The following lemma is a result from standard wavelet theory. See, for example, Daubechies (1992).

LEMMA 4. *Suppose the wavelet $\psi$ has compact support and is $r$-regular with $r > \alpha$. Then there exists a constant $C > 0$ such that for all $f \in F(\alpha, M)$ its wavelet coefficients satisfy*

$$|\theta_{j,k}| \leq C 2^{-j((1/2)+\alpha)} \qquad \text{for all } j \geq j_0 \text{ and } 1 \leq k \leq 2^j. \tag{46}$$

### 6.2. Proof of the main results.

PROOF OF THEOREM 1.  The proof of this theorem is divided into two parts. In the first part lower bounds are given and in the second upper bounds. For the lower bounds only the first two cases in the theorem need to be considered. Since the proofs of these two cases are similar a proof of case (i) is given in detail and then only the main changes needed for the proof of case (ii) are given.

*Lower bounds.*

CASE (i).  Let $g : \mathbb{R} \to \mathbb{R}$ be a function satisfying:

(i) $g(x) = \lambda > 0$ for $x \in [-1, 1]$ and $g$ is compactly supported in the interval $[-A, A]$;
(ii) $|g^{(k)}(x) - g^{(k)}(y)| \leq (M - M')|x - y|^{\alpha - k}, -\infty < x < y < \infty$ where $k$ is the greatest integer less than or equal to $\alpha$;
(iii) $\int_{-A}^{A} g^2(x)\, dx = 1$.

For sufficiently large $A$ such a function is easy to construct.

Set
$$\gamma_n = \left(\frac{n}{\log B_n}\right)^{\alpha/(1+2\alpha)} \quad \text{and} \quad \beta_n = \left(\frac{n}{\log B_n}\right)^{1/(1+2\alpha)}$$

and note that
$$\beta_n \gamma_n^2 = \frac{n}{\log B_n} \quad \text{and} \quad \beta_n^\alpha \gamma_n^{-1} = 1.$$

Let $f_{n,\theta}:[0,1] \to \mathbb{R}$ be defined by
$$f_{n,\theta}(x) = \theta \cdot \gamma_n^{-1} g(\beta_n(x - x_0)) + f_0(x) \quad \text{for } \theta = 0, 1.$$

It is simple to check that for $\theta = 0$ or $1$, $f_{n,\theta} \in F(\alpha, M)$ for all $n$. Note also that for sufficiently large $n$, say $n \geq N_0$,

(47) $$\rho_n = n \int_0^1 (f_{n,1} - f_{n,0})^2 = \log B_n.$$

Write $P_\theta^n$ for the probability measure associated with the process
$$Z_n^*(t) \equiv \int_0^t f_{n,\theta}(x)\,dx + \frac{1}{\sqrt{n}} B^*(t), \quad 0 \leq t \leq 1.$$

A sufficient statistic for the family of measures $\{P_0^n, P_1^n\}$ is then given by the log likelihood ratio $T_n = \ln \frac{dP_1^n}{dP_0^n}$, and for $n \geq N_0$,

$$\text{under } P_0^n, T_n \sim N\left(-\frac{\rho_n}{2}, \rho_n\right)$$

and

$$\text{under } P_1^n, T_n \sim N\left(\frac{\rho_n}{2}, \rho_n\right).$$

Now based on the Gaussian model (2) let $\hat{f}_n$ be an estimator of $f$. Decompose this estimator into components

(48) $$\hat{f}_n(x) = f_{n,0}(x) + \hat{\theta}(f_{n,1}(x) - f_{n,0}(x)) + \hat{h}_n(x),$$

where

(49) $$\int_{x_0 - c_n}^{x_0 + c_n} \hat{h}_n(x)(f_{n,1}(x) - f_{n,0}(x))\,dx = 0.$$

Hence, for $\theta = 0$ or $1$,

(50) $$\begin{aligned}\frac{1}{2c_n} \int_{x_0-c_n}^{x_0+c_n} (\hat{f}_n(x) - f_{n,\theta}(x))^2\,dx \\ \geq (\hat{\theta} - \theta)^2 \frac{1}{2c_n} \int_{x_0-c_n}^{x_0+c_n} (f_{n,1}(x) - f_{n,0}(x))^2\,dx \\ = (\hat{\theta} - \theta)^2 \beta_n^{-1} \gamma_n^{-2} \frac{1}{2c_n} \int_{-\beta_n c_n}^{\beta_n c_n} g^2(x)\,dx.\end{aligned}$$

It follows from the condition $\limsup_{n\to\infty} d_n \cdot (\log B_n)^{-1/(1+2\alpha)} = 0$ that for sufficiently large $n$, say $n \geq N_1$,
$$d_n \leq (\log B_n)^{1/(1+2\alpha)},$$
in which case $\beta_n c_n \leq 1$. So

$$(51) \quad \frac{1}{2c_n} \int_{x_0-c_n}^{x_0+c_n} (\hat{f}_n(x) - f_{n,\theta}(x))^2 \, dx \geq (\hat{\theta} - \theta)^2 \lambda^2 \left(\frac{\log B_n}{n}\right)^{2\alpha/(1+2\alpha)}.$$

If assumption (5) of the theorem holds, there exist a $C_1 < \infty$ and $N_2$ such that for all $n \geq N_2$,
$$R(\hat{f}_n, f_{n,0}; x_0, c_n) \leq C_1 n^{-2\alpha/(1+2\alpha)} B_n^{-1}.$$

Hence
$$E_{f_{n,0}}(\hat{\theta} - 0)^2 \leq C_1 \lambda^{-2} B_n^{-1} (\log B_n)^{-2\alpha/(1+2\alpha)}.$$

Since $T_n$ is sufficient for $\{P_0^n, P_1^n\}$ apply Theorem 1 of Brown and Low (1996) with $I = e^{\rho_n} = B_n$ for $n \geq N_0$.

Let $N = \max(N_0, N_1, N_2)$. Theorem 1, equation (2.4), of Brown and Low (1996) then yields for $n \geq N$

$$(52) \quad E_{f_{n,1}}(\hat{\theta} - 1)^2 \geq 1 - 2C_1^{1/2} \lambda^{-1} (\log B_n)^{-\alpha/(1+2\alpha)}.$$

Combining (51) and (52) yields (6).

CASE (ii). Let the function $g$ be constructed similarly as in the proof of case (i), except that in Condition (ii), $M - M'$ is replaced by $(M - M')(\frac{d}{A})^{(1/2)+\alpha}$. Set

$$\gamma_n = \left(\frac{d}{A}\right)^{1/2} \left(\frac{n}{\log B_n}\right)^{\alpha/(1+2\alpha)} \quad \text{and}$$

$$\beta_n = \frac{A}{d} \left(\frac{n}{\log B_n}\right)^{1/(1+2\alpha)}.$$

Since $\liminf_{n\to\infty} d_n \cdot (\log B_n)^{-1/(1+2\alpha)} > 0$, there exist constants $d > 0$ and $N > 0$ such that for all $n > N$,
$$d_n \geq d(\log B_n)^{1/(1+2\alpha)}.$$

Hence, for $n > N$,
$$\beta_n \gamma_n^2 = \frac{n}{\log B_n}, \qquad \beta_n^\alpha \gamma_n^{-1} = \left(\frac{A}{d}\right)^{(1/2)+\alpha} \quad \text{and} \quad \beta_n c_n \geq A.$$

In this case (50) yields
$$\frac{1}{2c_n} \int_{x_0-c_n}^{x_0+c_n} (\hat{f}_n(x) - f_{n,\theta}(x))^2 \, dx \geq \frac{1}{2}(\hat{\theta} - \theta)^2 n^{-2\alpha/(1+2\alpha)} \cdot \frac{\log B_n}{d_n}.$$

The remaining steps are the same as in the proof of Case (i) and hence are omitted.

We now turn to the proof of upper bounds, where the three cases need to be treated separately. Note, however, that in each case we may assume without loss of generality that $f_0 \equiv 0$ since we can always recenter the estimate at any given $f_0$. Let $\{\phi, \psi\}$ be a pair of compactly supported father and mother wavelets generating an orthonormal basis in $L^2[0,1]$ where the support lengths of $\phi$ and $\psi$ are denoted by $N_\phi$ and $N_\psi$, respectively. We assume that both $\phi$ and $\psi$ have $r > \alpha$ vanishing moments, $\int x^k \phi(x)\, dx = 0$ for $k = 1, \ldots, r$ and $\int x^k \psi(x)\, dx = 0$ for $k = 0, 1, \ldots, r$. For example, Coiflets of order greater than $\alpha$ have this property. See Daubechies (1992).

*Upper bounds.*

CASE (i). Let $j_n$ be the largest integer satisfying $2^{j_n} \leq (\frac{n}{\log B_n})^{1/(1+2\alpha)}$. For $j \geq 0$ and $1 \leq k \leq 2^j$ let $\phi_{j,k}(t) = 2^{j/2}\phi(2^j t - k)$. Then $x_0 \in \text{supp}(\phi_{j_n,k})$ for some $k$. Write

$$\tilde{y}_n \equiv 2^{j_n/2} \int \phi_{j_n,k}(t)\, dZ_n^*(t)$$
$$= 2^{j_n/2} \int f(t)\phi_{j_n,k}(t)\, dt + 2^{j_n/2} n^{-1/2} \int \phi_{j_n,k}(t)\, dW(t)$$
$$\equiv \bar{f} + z.$$

Here $z$ is a Gaussian random variable with mean 0 and variance $\sigma_n^2 = 2^{j_n} n^{-1}$, and $\bar{f}$ can be regarded as the "mean value" of $f$ on the support of $\phi_{j_n,k}$. Set

$$\delta_n = \text{sgn}(\tilde{y}_n)(|\tilde{y}_n| - \sigma_n (2 \log B_n)^{1/2})_+$$

and let $\hat{f}_n$ be an estimator of $f$ with

$$\hat{f}_n(x) \equiv \delta_n \qquad \text{for all } x \in [x_0 - c_n, x_0 + c_n].$$

We show below that $\hat{f}_n$ satisfies both $\hat{f}_n \in \Delta(f_0)$ and (7). First, it is easy to verify directly that

$$E(\delta_n - \bar{f})^2 \leq \min(2(\bar{f})^2, \sigma_n^2(1 + 2\log B_n)) + \sigma_n^2 B_n^{-1}$$

and hence

(53)
$$R(\hat{f}_n, f; x_0, c_n)$$
$$= \frac{1}{2c_n} E_f \int_{x_0 - c_n}^{x_0 + c_n} (\delta_n - f(x))^2\, dx$$
$$\leq \frac{1}{2c_n} \int_{x_0 - c_n}^{x_0 + c_n} (\bar{f} - f(x))^2\, dx$$
$$\quad + \min(2(\bar{f})^2, \sigma_n^2(1 + 2\log B_n)) + \sigma_n^2 B_n^{-1}.$$

Now for $f_0 \equiv 0$ the zero function on $[0,1]$ the first two terms in (53) are both 0. Hence

$$R(\hat{f}_n, f_0; x_0, c_n) \leq n^{-2\alpha/(1+2\alpha)} B_n^{-1} (\log B_n)^{-1/(1+2\alpha)}$$

and it follows that $\hat{f}_n \in \Delta(f_0)$. It follows from the vanishing moments property of $\phi$ that for all $f \in F(\alpha, M)$ and for all $x \in [x_0 - c_n, x_0 + c_n]$,

(54) $$|f(x) - \bar{f}| \leq C(M, \phi) 2^{-\alpha j_n},$$

where $C(M, \phi)$ is a constant depending on $M$ and $\phi$ only. Now (7) follows by applying (54) to (53):

$$R(\hat{f}_n, f; x_0, c_n) \leq C \left( \frac{\log B_n}{n} \right)^{2\alpha/(1+2\alpha)} (1 + o(1)).$$

CASE (ii). In the second case, a wavelet procedure based on block thresholding is used. Let $J_1$ and $J_2$ be the largest integers satisfying

$$2^{J_1} \leq d_n^{-1} n^{1/(1+2\alpha)} \quad \text{and} \quad 2^{J_2} \leq d_n^{-1} \log B_n n^{1/(1+2\alpha)}.$$

Let

$$H_j = \{(j,k) : \operatorname{supp}(\psi_{j,k}) \cap [x_0 - c_n, x_0 + c_n] \neq \varnothing\} \quad \text{and} \quad H_* = \bigcup_{J_1 \leq j \leq J_2} H_j.$$

Then it is easy to check that for $j \geq J_1$ the cardinality of the index sets $H_j$ is of order $2^j c_n$ and so $L_n \equiv \operatorname{Card}(H_*) = b_n \log B_n$ with $b_* \leq b_n \leq b^*$ for some positive constants $b_*$ and $b^*$. Denote by $S^2$ the sum of all the squared empirical wavelet coefficients $y_{j,k}$ with indices in $H_*$. Applying a block thresholding rule to the coefficients,

$$\hat{\theta}_{j,k} = \left( 1 - \frac{\lambda L_n n^{-1}}{S^2} \right)_+ y_{j,k} \quad \text{for all } (j,k) \in H_*.$$

Then it follows from Lemma 2 that

(55)
$$\sum_{(j,k) \in H_*} E(\hat{\theta}_{j,k} - \theta_{j,k})^2$$
$$\leq \min \left( \sum_{(j,k) \in H_*} \theta_{j,k}^2, \lambda L_n n^{-1} \right) + 2n^{-1} e^{-(1/2)(\lambda - \log \lambda - 1) L_n}.$$

Let the thresholding constant $\lambda$ be chosen such that

$$\tfrac{1}{2}(\lambda - \log \lambda - 1) b_* = 1.$$

Then the second term in the right-hand side of (55) is bounded from above by $2n^{-1}B_n^{-1}$. Applying Lemma 1, we have

$$R(\hat{f}_n, f; x_0, c_n)$$

(56)
$$\leq \frac{1}{2c_n} \sum_{(j,k)\in H_*} E(\hat{\theta}_{j,k} - \theta_{j,k})^2 + \frac{1}{2c_n} \sum_{j>J_2} \sum_{(j,k)\in H_j} \theta_{j,k}^2$$

$$\leq \frac{1}{2c_n} \min\left(\sum_{(j,k)\in H_*} \theta_{j,k}^2, \lambda L_n n^{-1}\right)$$

$$+ n^{-2/(1+2\alpha)} B_n^{-1} (\log B_n)^{-1/(1+2\alpha)} + \frac{1}{2c_n} \sum_{j>J_2} \sum_{(j,k)\in H_j} \theta_{j,k}^2.$$

For $f = f_0 \equiv 0$, the first and the third terms in (56) are both 0, hence

$$R(\hat{f}_n, f_0; x_0, c_n) \leq n^{-2\alpha/(1+2\alpha)} B_n^{-1} (\log B_n)^{-1/(1+2\alpha)}$$

and so $\hat{f}_n \in \Delta(f_0)$. For $f \in F(\alpha, M)$, it follows from Lemma 4 that

(57)
$$|\theta_{j,k}| \leq C 2^{-j((1/2)+\alpha)}$$

with the constant $C$ not depending on $f$. Hence

(58)
$$\frac{1}{2c_n} \sum_{j>J_2} \sum_{(j,k)\in H_j} \theta_{j,k}^2 \leq \frac{1}{2c_n} \sum_{j>J_2} C 2^j c_n 2^{-j(1+2\alpha)}$$

$$= C n^{-2\alpha/(1+2\alpha)} d_n^{-1} \log B_n.$$

Now (9) follows from (56) and (58).

CASE (iii). Finally, we turn to the third case where we will use the same notation as in Case (ii). Let $J_2$ and $J_3$ be the largest integers satisfying $2^{J_2} \leq d_n^{-1} \log B_n n^{1/(1+2\alpha)}$ and $2^{J_3} \leq n^{1/(1+2\alpha)}$, respectively. (If $d_n < \log B_n$, choose $J_3 = J_2$.) Denote by $L_j$ the cardinality of the index sets $H_j$. Then, for $j \geq J_2$, there exist positive constants $b_*$ and $b^*$ such that $b_* 2^j c_n \leq L_j \leq b^* 2^j c_n$. Denote by $S_j^2$ the sum of all the squared empirical wavelet coefficients $y_{j,k}$ at level $j$ with $(j,k) \in H_j$. Applying a block thresholding rule to the coefficients level by level,

$$\hat{\theta}_{j,k} = \left(1 - \frac{\lambda L_j n^{-1}}{S_j^2}\right)_+ y_{j,k} \quad \text{for all } J_2 \leq j \leq J_3 \text{ and } (j,k) \in H_j.$$

Then again it follows from Lemma 2 that

(59)
$$\sum_{(j,k)\in H_j} E(\hat{\theta}_{j,k} - \theta_{j,k})^2$$

$$\leq \min\left(\sum_{(j,k)\in H_j} \theta_{j,k}^2, \lambda L_j n^{-1}\right) + 2n^{-1} e^{-(1/2)(\lambda - \log\lambda - 1)L_j}.$$

Write $L_{J_2} = b_n \log B_n$ with $b_* \leq b_n \leq b^*$. We choose the thresholding constant $\lambda$ such that

$$\tfrac{1}{2}(\lambda - \log \lambda - 1)b_* = 1.$$

Then the second term on the right-hand side of (59) is bounded from above by $2n^{-1}B_n^{-1}$ for $j = J_2$ and

$$(60) \qquad \sum_{j=J_2}^{J_3} 4n^{-1}e^{-(1/2)(\lambda - \log \lambda - 1)L_j} \leq 4n^{-1}B_n^{-1}.$$

Lemma 1 yields

$$\begin{aligned}
R(\hat{f}_n, f; x_0, c_n) \\
\leq \frac{1}{2c_n} \sum_{j=J_2}^{J_3} \sum_{(j,k) \in H_j} E(\hat{\theta}_{j,k} - \theta_{j,k})^2 + \frac{1}{2c_n} \sum_{j>J_3} \sum_{(j,k) \in H_j} \theta_{j,k}^2 \\
\leq \frac{1}{2c_n} \sum_{j=J_2}^{J_3} \min\left(\sum_{(j,k) \in H_j} \theta_{j,k}^2, \lambda L_j n^{-1}\right) \\
+ 2n^{-2\alpha/(1+2\alpha)} B_n^{-1} d_n^{-1} + \frac{1}{2c_n} \sum_{j>J_3} \sum_{(j,k) \in H_j} \theta_{j,k}^2.
\end{aligned}$$

(61)

Once again for $f = f_0 \equiv 0$, the first and the third terms in (61) are both 0, hence

$$R(\hat{f}_n, f_0; x_0, c_n) \leq 2n^{-2\alpha/(1+2\alpha)} B_n^{-1} d_n^{-1}$$

and so $\hat{f}_n \in \Delta(f_0)$. The coefficient bound (57) yields

$$(62) \quad \frac{1}{2c_n} \sum_{j>J_3} \sum_{(j,k) \in H_j} \theta_{j,k}^2 \leq \frac{1}{2c_n} \sum_{j>J_3} \sum_{(j,k) \in H_j} Cc_n 2^{-2\alpha j} = Cn^{-2\alpha/(1+2\alpha)}.$$

Now (10) follows from (61) and (62). □

PROOF OF THEOREM 2. As in the proof of Theorem 1 it suffices to consider $f_0 \equiv 0$. Expand the Gaussian process (2) in an orthonormal wavelet basis as in Section 3.1. Suppose the wavelet $\psi$ is chosen to be $r$-regular with $r > \alpha$. Let $J'$ be the largest integer satisfying $2^{J'} < n^{1/(1+2\alpha)}$. Then the total number $L'$ of wavelet coefficients up to (and including) the level $J'$ is less than $2n^{1/(1+2\alpha)}$ and larger than or equal to $n^{1/(1+2\alpha)}$. Group all the empirical wavelet coefficients $\tilde{y}_{j_0,k}$ and $y_{j,k}$ up to the level $J'$ into a single block and apply a James–Stein type rule to the coefficients. More specifically, denote the sum of the squared empirical coefficients up to the level $J'$ by

$$S^2 = \sum_{k=1}^{2^{j_0}} \tilde{y}_{j_0,k}^2 + \sum_{j=j_0}^{J'} \sum_{k=1}^{2^j} y_{j,k}^2$$

and define the estimator of the wavelet coefficients by

$$\hat{\xi}_{j_0,k} = \left(1 - \frac{\lambda L' n^{-1}}{S^2}\right)_+ \tilde{y}_{j_0,k} \quad \text{for } 1 \leq k \leq 2^{j_0},$$

(63) $$\hat{\theta}_{j,k} = \left(1 - \frac{\lambda L' n^{-1}}{S^2}\right)_+ y_{j,k} \quad \text{for } j \leq J', 1 \leq k \leq 2^j,$$

$$\hat{\theta}_{j,k} = 0 \quad \text{otherwise},$$

where $\lambda$ is a constant satisfying $\lambda - \log \lambda - 1 = 2D$. The corresponding estimator $\hat{f}_n$ of $f$ is the wavelet series with $\hat{\xi}_{j_0,k}$ and $\hat{\theta}_{j,k}$ as coefficients:

(64) $$\hat{f}_n(x) = \sum_{k=1}^{2^{j_0}} \hat{\xi}_{j_0,k} \phi_{j_0,k}(x) + \sum_{j=j_0}^{\infty} \sum_{k=1}^{2^j} \hat{\theta}_{j,k} \psi_{j,k}(x).$$

It follows from (42) in Lemma 2 that

$$\sum_k E(\hat{\xi}_{j_0,k} - \xi_{j_0,k})^2 + \sum_{j=j_0}^{J'} \sum_k E(\hat{\theta}_{j,k} - \theta_{j,k})^2$$

$$\leq \min\left(\sum_k \xi_{j_0,k}^2 + \sum_{j=j_0}^{J'} \sum_k \theta_{j,k}^2, \lambda L' n^{-1}\right) + 2n^{-1} e^{-(1/2)(\lambda - \log \lambda - 1)L'}$$

$$\leq \min\left(\sum_k \xi_{j_0,k}^2 + \sum_{j=j_0}^{J'} \sum_k \theta_{j,k}^2, 2\lambda n^{-2\alpha/(1+2\alpha)}\right) + 2n^{-1} e^{-Dn^{1/(1+2\alpha)}}.$$

Hence

$$E_f \|\hat{f}_n - f\|_2^2$$

$$= \sum_k E(\hat{\xi}_{j_0,k} - \xi_{j_0,k})^2 + \sum_{j=j_0}^{J'} \sum_k E(\hat{\theta}_{j,k} - \theta_{j,k})^2 + \sum_{j=J'+1}^{\infty} \sum_k \theta_{j,k}^2$$

(65) $$\leq \min\left(\sum_k \xi_{j_0,k}^2 + \sum_{j=j_0}^{J'} \sum_k \theta_{j,k}^2, 2\lambda n^{-2\alpha/(1+2\alpha)}\right)$$

$$+ 2n^{-1} e^{-Dn^{1/(1+2\alpha)}} + \sum_{j=J'+1}^{\infty} \sum_k \theta_{j,k}^2.$$

Now for $f_0 \equiv 0$, all $\xi_{j_0,k} = 0$ and all $\theta_{j,k} = 0$ so

$$E_{f_0} \|\hat{f}_n - f_0\|_2^2 \leq 2n^{-1} e^{-Dn^{1/(1+2\alpha)}}.$$

Thus, with $B_n = n e^{Dn^{1/(1+2\alpha)}}$,

$$\limsup_{n \to \infty} B_n E_{f_0} \|\hat{f}_n - f_0\|_2^2 < \infty \quad \text{and} \quad \lim_{n \to \infty} n^{1/(1+2\alpha)} (\log B_n)^{-1} = D.$$

On the other hand, the estimator attains the optimal rate uniformly over $F(\alpha, M)$. This can be seen easily from (46) and (65):

$$\sup_{f \in F(\alpha,M)} E_f \|\hat{f}_n - f\|_2^2$$

$$\leq 2\lambda n^{-2\alpha/(1+2\alpha)} + 2n^{-1} e^{-Dn^{1/(1+2\alpha)}} + \sum_{j=J'+1}^{\infty} \sum_{k=1}^{2^j} C^2 2^{-j(1+2\alpha)}$$

$$\leq 2(\lambda + C^2) n^{-2\alpha/(1+2\alpha)} (1 + o(1)). \qquad \Box$$

PROOF OF THEOREM 3. We assume $J^* < J$ in the following proof. In the special case of $J^* \geq J$ the estimator is the BlockJS estimator. The proof for this case follows from that of Theorem 4. Denote by $I_n(x) = I(x \in [x_0 - c_n, x_0 + c_n])$. Then

$$R(\hat{f}_n, f; x_0, c_n)$$

$$= E\left\{\frac{1}{2c_n} \int_0^1 \left[\sum_k (\tilde{y}_{j_0,k} - \xi_{j_0,k})\phi_{j_0,k}(x)\right.\right.$$

$$\left.\left. + \sum_{j=j_0}^{\infty} \sum_k (\hat{\theta}_{j,k} - \theta_{j,k})\psi_{j,k}(x)\right]^2 I_n(x)\, dx\right\}$$

$$\leq 2^{j_0} E\left\{\frac{1}{c_n} \int_0^1 \sum_k (\tilde{y}_{j_0,k} - \xi_{j_0,k})^2 \phi_{j_0,k}^2(x) I_n(x)\, dx\right\}$$

$$+ E\left\{\frac{1}{c_n} \int_0^1 \left[\sum_{j=j_0}^{\infty} \sum_k (\hat{\theta}_{j,k} - \theta_{j,k})\psi_{j,k}(x)\right]^2 I_n(x)\, dx\right\}$$

$$\leq Cn^{-1} + E\left\{\frac{1}{c_n} \int_0^1 \left[\sum_{j=j_0}^{\infty} \sum_k (\hat{\theta}_{j,k} - \theta_{j,k})\psi_{j,k}(x)\right]^2 I_n(x)\, dx\right\}.$$

Hence

$$R(\hat{f}_n, f; x_0, c_n)$$

$$\leq Cn^{-1} + 2\|\psi\|_\infty E\left(\sum_{j=j_0}^{J_*} \sum_{(j,k) \in H_j} 2^{j/2} |\hat{\theta}_{j,k} - \theta_{j,k}|\right)^2$$

(66) $$+ E\left\{\frac{2}{c_n} \int_0^1 \left(\sum_{j>J_*} \sum_k (\hat{\theta}_{j,k} - \theta_{j,k})\psi_{j,k}(x) I_n(x)\right)^2 dx\right\}$$

$$\leq Cn^{-1} + 2\|\psi\|_\infty \left(\sum_{j=j_0}^{J_*} \sum_{(j,k)\in H_j} 2^{j/2}(E(\hat{\theta}_{j,k} - \theta_{j,k})^2)^{1/2}\right)^2$$

$$+ \frac{2}{c_n} \sum_{j>J_*} \sum_{(j,k)\in H_j} E(\hat{\theta}_{j,k} - \theta_{j,k})^2.$$

The last inequality follows from Lemma 1 and the elementary inequality

$$E\left(\sum_{i=1}^n X_i\right)^2 \leq \left(\sum_{i=1}^n (EX_i^2)^{1/2}\right)^2.$$

We now consider the three cases separately. The main tool is the risk bounds (43) and (44) given in Lemma 2. Note that with $\sigma^2 = n^{-1}$ and $L = \log n$ the second term on the right-hand side of (43) and (44) is $2\lambda n^{-2}$, which is negligible in the following risk calculations, and we will absorb this term into the first term, $Cn^{-1}$, in the calculations below. Note that

$$\sum_{j=J_*}^{J^*-1} \sum_{(j,k)\in H_j} E(\hat{\theta}_{j,k} - \theta_{j,k})^2 \leq C \min\left((\log n)n^{-1}, \sum_{(j,k)\in H_j} \theta_{j,k}^2\right) + O(n^{-2}).$$

In case (i), let $J_0$ be the smallest integer satisfying $2^{J_0} \geq (\frac{n}{\log n})^{1/(1+2\alpha)}$. Then $J_0 < J_*$. It follows from Lemmas 2 and 4 that

$$R(\hat{f}_n, f; x_0, c_n)$$

$$\leq C\left(\sum_{j=j_0}^{J_0-1} 2^{j/2}(\log n)^{1/2}n^{-(1/2)} + \sum_{j=J_0}^{J_*-1} 2^{j/2}2^{-j((1/2)+\alpha)}\right)^2$$

$$+ Cc_n^{-1} \sum_{j=J_*}^{\infty} 2^j c_n 2^{-j(1+2\alpha)}$$

$$\leq C\left(\frac{\log n}{n}\right)^{2\alpha/(1+2\alpha)}.$$

In case (ii), Lemmas 2 and 4 yield that

$$R(\hat{f}_n, f; x_0, c_n) \leq C\left(\sum_{j=j_0}^{J_*-1} 2^{j/2}(\log n)^{1/2}n^{-(1/2)}\right)^2$$

$$+ Cc_n^{-1}(\log n)n^{-1} + Cc_n^{-1}\sum_{j=J^*}^{\infty} 2^j c_n 2^{-j(1+2\alpha)}$$

$$\leq C\frac{\log n}{d_n} n^{-2\alpha/(1+2\alpha)}.$$

In case (iii) let $J_1$ be the smallest integer satisfying $2^{J_1} \geq n^{1/(1+2\alpha)}$. We have

$$R(\hat{f}_n, f; x_0, c_n) \leq C \left( \sum_{j=j_0}^{J_*-1} 2^{j/2}(\log n)^{1/2} n^{-(1/2)} \right)^2 + Cc_n^{-1}(\log n) n^{-1}$$

$$+ Cc_n^{-1} \sum_{j=J_*}^{J_1-1} \frac{2^j c_n}{\log n}(\log n) n^{-1} + Cc_n^{-1} \sum_{j=J_1}^{\infty} 2^j c_n 2^{-j(1+2\alpha)}$$

$$\leq C n^{-2\alpha/(1+2\alpha)}. \qquad \square$$

PROOF OF THEOREM 4. Let $\hat{f}_n(x)$ be the BlockJS estimator given in (33). Denote $I_n(x) = I(x \in [x_0 - c_n, x_0 + c_n])$. Similarly as in (66), in the proof of Theorem 3, for any $T \geq j_0$,

$$R(\hat{f}_n, f; x_0, c_n) \leq Cn^{-1} + 2\|\psi\|_\infty \left( \sum_{j \leq T} \sum_{(j,k) \in H_j} 2^{j/2} (E(\hat{\theta}_{j,k} - \theta_{j,k})^2)^{1/2} \right)^2$$

(67)

$$+ \frac{2}{c_n} \sum_{j > T} \sum_{(j,k) \in H_j} E(\hat{\theta}_{j,k} - \theta_{j,k})^2.$$

Denote by $J_i$, $i = 0, 1, 2, 3, 4$, the smallest integers satisfying

$$2^{J_0} \geq \left( \frac{n}{\log n} \right)^{1/(1+2\alpha)}, \qquad 2^{J_1} \geq \frac{n^{1/(1+2\alpha)}}{d_n},$$

$$2^{J_2} \geq n^{1/(1+2\alpha)} \left( \frac{d_n}{\log n} \right)^{1/(2\alpha)},$$

$$2^{J_3} \geq \frac{n^{1/(1+2\alpha)} \log n}{d_n}, \qquad 2^{J_4} \geq n^{1/(1+2\alpha)}.$$

Then for all $j \leq J_1$,

$$\text{Card}(H_j) \leq \begin{cases} N_\psi, & \text{if } j < J_1, \\ N_\psi 2^j c_n, & \text{if } j \geq J_1. \end{cases}$$

Note that for all levels $j \leq J_3$, the coefficients of wavelet basis functions $\psi_{j,k}$ whose support has nonempty intersection with the interval $[x_0 - c_n, x_0 + c_n]$ are in at most $N_\psi + 1$ blocks because the number of such coefficients is less than $N_\psi \log n$.

We will consider the three cases separately. Again, with $\sigma^2 = n^{-1}$ and $L = \log n$, the second term on the right-hand side of (44) and (43) is $2\lambda n^{-2}$, which is negligible and thus will be absorbed into the $Cn^{-1}$ term in the calculations below.

(i) Choose $T = J_1$ in (67). In this case $J_0 < J_1$. It then follows from Lemmas 2 and 4 that

$R(\hat{f}_n, f; x_0, c_n)$

$\leq Cn^{-1} + C\left(\sum_{j=j_0}^{J_0-1} 2^{j/2}(\log n)^{1/2}n^{-(1/2)} + \sum_{j=J_0}^{J_1-1} 2^{j/2}2^{-j((1/2)+\alpha)}\right)^2$

$+ c_n^{-1}\sum_{j=J_1}^{\infty} 2^j c_n 2^{-j(1+2\alpha)}$

$\leq C\left(\frac{\log n}{n}\right)^{2\alpha/(1+2\alpha)}.$

(ii) Choose $T = J_1$ in (67). Lemmas 2 and 4 yield that

$R(\hat{f}_n, f; x_0, c_n)$

$\leq Cn^{-1} + C\left(\sum_{j=j_0}^{J_1-1} 2^{j/2}(\log n)^{1/2}n^{-(1/2)}\right)^2 + Cc_n^{-1}\sum_{j=J_1}^{J_2-1}(\log n)n^{-1}$

$+ Cc_n^{-1}\sum_{j=J_2}^{\infty} 2^j c_n 2^{-j(1+2\alpha)}$

$\leq Cn^{-1} + C\frac{\log n}{d_n}n^{-2\alpha/(1+2\alpha)}$

$+ C(J_2 - J_1)\frac{\log n}{d_n}n^{-2\alpha/(1+2\alpha)} + C\frac{\log n}{d_n}n^{-2\alpha/(1+2\alpha)}$

$= C\frac{(\log n)(\log\log n)}{d_n}n^{-2\alpha/(1+2\alpha)}(1 + o(1)).$

(iii) Choose $T = J_3$ in (67). In this case we have

$R(\hat{f}_n, f; x_0, c_n)$

$\leq Cn^{-1} + C\left(\sum_{j=j_0}^{J_1-1} 2^{j/2}(\log n)^{1/2}n^{-1/2}\right)^2 + Cc_n^{-1}\sum_{j=J_1}^{J_3-1}(\log n)n^{-1}$

$+ Cc_n^{-1}\sum_{j=J_3}^{J_4-1}\frac{2^j c_n}{\log n}(\log n)n^{-1} + Cc_n^{-1}\sum_{j=J_4}^{\infty} 2^j c_n 2^{-j(1+2\alpha)}$

$\leq Cn^{-1} + C\frac{\log n}{d_n}n^{-2\alpha/(1+2\alpha)}$

$+ C(J_3 - J_1)\frac{\log n}{d_n}n^{-2\alpha/(1+2\alpha)} + Cn^{-2\alpha/(1+2\alpha)} + Cn^{-2\alpha/(1+2\alpha)}$

$$= Cn^{-2\alpha/(1+2\alpha)}(1+o(1)). \qquad \square$$

**Acknowledgments.** We thank the referee and an Associate Editor for thorough and useful comments which have helped to improve the presentation of the paper.

Department of Statistics  
The Wharton School  
University of Pennsylvania  
Philadelphia, Pennsylvania 19104-6340  
USA  
e-mail: tcai@wharton.upenn.edu